\documentclass{amsart}

\newtheorem{theorem}{Theorem}[section]
\newtheorem{lemma}[theorem]{Lemma}
\newtheorem{corollary}[theorem]{Corollary}
\newtheorem{prop}[theorem]{Proposition}
\theoremstyle{definition}
\newtheorem{definition}[theorem]{Definition}

\theoremstyle{remark}
\newtheorem{remark}[theorem]{Remark}
\numberwithin{equation}{section}

\def\G{{\bf G}}
\def\Z{{\bf Z}}
\def\H{{\bf H}}

\def\N{{\bf N}}

\def\Spec{\mathop{\rm Spec}\nolimits}
\def\Nm{\mathop{\rm Nm}\nolimits}

\def\Ind{\mathop{\rm Ind}\nolimits}
\def\res{\mathop{\rm res}\nolimits}
\def\red{\mathop{\rm red}\nolimits}
\def\ind{\mathop{\rm ind}\nolimits}
\def\Im{\mathop{\rm Im}\nolimits}
\def\Aut{\mathop{\rm Aut}\nolimits}

\def\Def{\mathop{\rm Def}\nolimits}

\def\Set{\mathop{\rm Set}\nolimits}

\begin{document}
\title{Induction and restriction in formal deformation of coverings}
\author{Jos\'e Bertin and Ariane M\'ezard}
\email{jose.bertin@ujf-grenoble.fr, ariane.mezard@math.u-psud.fr}
\maketitle
\noindent
{\bf Abstract}: 
Let $X/S$ be a semistable curve with an action of a finite 
group $G$ and let $H$ be a normal subgroup of $G$. We present a new 
condition under which for any base change $T\rightarrow S$, 
$(X/G)\times_S T$ is isomorphic to $(X\times_S T)/G$. This allows us
 to define induction and restriction morphisms between 
 the $G$-equivariant deformation functor of $X$ and the
 $G/H$-equivariant (resp. $H$-equivariant) deformation 
functor of $X/H$ (resp. $X$).\\

\noindent
{\bf R\'esum\'e} : Soit $X/S$ une courbe semi-stable munie de l'action  
d'un groupe fini $G$. Soit $H$ un sous-groupe normal de $G$. Nous proposons une 
nouvelle condition sous laquelle $(X/G)\times_S T$ est isomorphe \`a 
$(X\times_S T)/G$ quel que soit le changement de base $T\rightarrow S$. Ce
r\'esultat nous permet de d\'efinir des morphismes d'induction et de 
restriction entre le foncteur de d\'eformations 
$G$-\'equivariantes de $X$ et   
le foncteur de d\'eformations $G/H$-\'equivariantes (resp. 
$H$-\'equivariantes) de
$X/H$ (resp. $X$).

\section{Introduction}
Let $S = \Spec R$ be a base scheme, where $R$ is a commutative local 
noetherian ring. Let $X \to S$ be a semistable curve with a  
faithfull action of a
finite group $G$, and let  $H\triangleleft\, G$ 
be a normal subgroup. There is an induced action of
$G/H$ on the quotient $X/H$
which makes the quotient map $X \to X/H$ equivariant with respect to
$G \to G/H$. 
We wish to compare the  $G$-equivariant deformation functor 
$\Def_G(X)$ of $X$ (in the sense of \cite{BeMe})
with the
corresponding $G/H$-equivariant deformation
functor $\Def_{G/H}(X/H)$ of $X/H$. This amounts to establish  
the existence of a morphism  between both deformation functors. 
This will be realized provide that base change and passage to
 quotient are commuting operations. This means that for a base change 
$T \to S$  
the natural morphism 
\begin{equation}
\label{morphi} 
(X/G)\times_S T \longrightarrow (X\times_S T)/G 
\end{equation} is an isomorphism. 
After (\cite{KaMa} Theorem 7.1.3), we know (\ref{morphi}) 
 universally is an isomorphism if one of the following two 
conditions is satisfied:
\hfill\break $\bullet$ $G$ acts freely 
\hfill\break $\bullet$ the order of $G$ is invertible in the
 structural sheaf of $S$.
\hfill\break 
Our first result is a new condition under which  (\ref{morphi})  
universally is an isomorphism:
\begin{theorem}
\label{th11}
 Let $X \rightarrow S$ be a semistable curve, and $G$ be a finite group of $S$-automorphisms of $X/S$. If the action of $G$ is free on an open dense set on any geometric fiber, then (\ref{morphi}) is universally an isomorphism. 
\end{theorem}  
In addition we provide a local version (Proposition \ref{localprop})  
of theorem \ref{th11} for $\G$ a finite flat $S$-group scheme acting on
 a smooth affine $S$-curve. The reason for considering in the 
local case a finite flat  group scheme is the following: under 
a stable specialisation in characteristic $p>0$ dividing the order
 of a Galois group $G$ of the  covering, the ordinary group $G$  may 
specialize, at least locally, in an infinitesimal group scheme. 
This suggests that we must take into account in the general 
problem of deformation-specialisation of coverings,  not only
 constant groups, but also infinitesimal groups.\\

What is important 
for application of theorem \ref{th11} to deformations, is that there is
  no condition on the residue characteristic of $S$, relatively on
 the order of $G$.   Under theorem \ref{th11}'s hypotheses, 
we construct  
 induction and restriction morphisms 
\begin{equation} 
\label{morDef} 
\Def_{H}(X)\stackrel \res\longleftarrow \Def_G (X) \stackrel 
\ind \longrightarrow \Def_{G/H} (X/H) 
\end{equation}

The outline of this paper is as follows: first  we present a version
 of the Kleiman-Lonsted algebra (\cite{KlLo}) in a slightly generalized form,
 i.e. to  a finite flat $S$-group  scheme admissibly acting on $X$ 
(\S \ref{paraKL}). 
By using the Kleiman-Lonsted algebra we prove 
 a base change criterion (\S \ref{critere}) and  a local version
 of \ref{th11} (\S \ref{localversion}).
 After having worked out  the double points case we
 prove theorem \ref{th11} (\S \ref{semistable}).  
 At last we define the  induction-restriction morphisms between 
the deformation functors and we compute the differentials of these morphisms
for local deformations functors
(\S \ref{dersection}).

 \section{Kleiman-Lonsted's algebra} 
\label{paraKL} 

In this section we briefly present  the basic 
properties of the Kleiman-Lonsted algebra ([KlLo]) associated 
to a finite flat group scheme action. This algebra  is an approximation 
of the algebra of invariants (\S \ref{subsection31}).  
But in the case of an infinitesimal group scheme acting on a singular curve, this algebra is more narrow.
We prove
some properties of reduction of this algebra 
(\S \ref{propriete}). Note that  most  of the proofs 
given below remain equally valid for an action of a finite 
flat groupo\"\i d (\cite{DeGa}), but due to lack of applications,
 we shall limit ourselves in the statement, to group  schemes actions.

\subsection{Definition and elementary properties} 
\label{subsection31} 
Let $R$ be a commutative noetherian ring.
 Let $\G$ be   finite  flat $R$-group scheme. We assume that the
 rank of $R[\G]$ over $R$ is constant and equals  the order $n = 
|\G|$ of $\G$. We consider an action of  $\G$ on a finite type 
$R$-algebra $A$.  The coaction is denoted by $\mu^*: A 
\rightarrow R[\G]\otimes_R A$. The norm mapping $ N_{\G}:
 A \rightarrow A$ is defined as usual (\cite{Mu},\cite {DeGa})
 by $N(a) = \Nm(\mu^* (a))$, where $\Nm$ is the norm of the 
$A$-algebra $R[\G]\otimes_R A$, the structural morphism being 
$a \to 1\otimes a$. The map $N_\G$ extends canonically to  
$N_\G:A[T]\rightarrow A[T]$. The characteristic polynomial 
of $a\in A$ is  $\chi _a (T) = N_\G(T - a)$. The $\G$-symmetric
 functions $\sigma_i^{\G}=\sigma_i$, $1\leq i\leq n$ 
 are defined by
 $$\chi_a (T) = T^n - \sigma_1(a)T^{n-1} + \cdots
 + (-1)^n \sigma_n(a)$$ 
After Cayley-Hamilton  theorem $\chi_a (a) = 0$ and 
$N(a) = \sigma_n (a)$. Similar definitions hold for an action 
of a finite flat groupo\"\i d (\cite{DeGa},\cite{Gi}).
 \begin{definition}
 The Kleiman-Lonsted algebra is defined as the $R$-algebra
 $$\Sigma_R^\G (A) = R[ \cup_{1\leq i\leq n} \sigma_i(A)]$$ 
\end{definition}

After \cite{Mu} \S12 Theorem 1, the norm maps $A$ into 
$A^\G=\{a\in  A,\mu^*(a)=1\otimes a\}$. The  $\sigma_i (a)$'s 
also belong to $A^\G$ and $\Sigma_R^\G (A) \subset A^\G$.
 We point out that as in the case of an action of a constant 
group (\cite{KlLo}), the algebras $\Sigma_R^\G(A)$ and $A^\G$
 are  closely related. We list  below some  basic facts about the
 Kleiman-Lonsted algebra. First, the $R$-algebra $\Sigma_R^\G(A)$ 
is of finite type over $R$. Second the $R$-algebras $A$ and $A^\G$
 are finite modules over $\Sigma_R^\G (A)$. Then the formation of $A^\G$
and $\Sigma_R^\G(A)$ commutes with flat base change: 
\begin{lemma}
 \label{propref} 
Let $R'$ be a $R$-algebra. Denote
 $\phi:\Sigma_R^{\G}(A)\otimes _RR'\rightarrow
 \Sigma_{R'}^{\G\otimes R'}(A\otimes_RR')$ and 
$\psi:A^\G\otimes_RR'\rightarrow (A\otimes_R R')^{\G\otimes R'}$.
Assume $\G$ is defined over an integral base ring.\\ 
 i. The morphism  $\phi$ is surjective .\\ 
 ii. If $R'$ is a flat $R$-algebra, $\psi$ and $\phi$ are bijective maps.\\ 
iii. if $R'$ is a faithfully flat $R$-algebra,
 or if the order of $\G$ is invertible, 
$\Sigma_R^\G(A)=A^\G$ if and only if 
$\Sigma_{R'}^{\G\otimes R'}(A\otimes_RR')=(A\otimes_RR')^{\G\otimes R'}$. 
\end{lemma} 
\begin{proof} 
The proof of assertions ii. and iii. in  ([KlLo] prop 4.4) works without  
any change provided i. holds.
  
First assume  $G$ to be a constant group. 
The proof in this case follows 
 from the  
following known fact
(see Appendix).
Let $n=|G|$ and $s_j(x)$, $1\leq j\leq n$ denote the usual symmetric functions of order $j$ in $n$ variables $x = (x_1,\ldots,x_n)$:
$$s_j(x)=\sum_{1\leq i_1<\cdots<i_j\leq n}x_{i_1}\cdots x_{i_j}$$
Let $x^{(\ell)}$, $1\leq \ell\leq q$ be $q$ sets of $n$ variables.
The quantity  $s_j(x^{(1)}+\cdots+x^{(q)})$ may be 
expressed as a polynomial with 
integer coefficients
in the symmetric functions $s_m, m\leq j$, where the arguments in the 
$s_m$'s are monomonials in $x^{(\ell)}$, $1\leq \ell\leq q$.

Denote by $(g_i)_{1\leq i\leq n}$ the elements of $G$. Since 
$\sigma_k(a)=s_k(g_1a,\ldots,g_na)$, $a\in A$, $1\leq k\leq n$,
then 
$$\forall  1\leq k\leq n, a\in A, r'\in R'\sigma_k(a\otimes r')=\sigma_k(a)\otimes r'\in \Sigma_{R'}^G(A\otimes_R
 R') $$
hence $\forall a'=\sum_{i=1}^ma_i\otimes r'_i\in A\otimes_RR',\;\;\; 
\sigma_k(a')\in\Sigma_{R'}^G(A\otimes_R
 R')$
and $\phi$ is surjective.\\

Second assume that 
$\G$ is etale over $R$. Perform an etale surjective (and thus faithfully flat)
base change and use the compatibility of the norm map to deduce the assertion
i. in an overing of $A$.

Then assume that 
$\G$ is not etale and $R=k$ is a perfect field of characteristic 
$p>0$. Denote by $\mu^*:A\rightarrow k[\G]\otimes_k A$ the coaction.
 The $k$-group scheme  $\G$ is the semi-direct product of its etale part 
$\G_{\red}$
 and its purely infinitesimal connected component
 $\G_0$ (\cite {Ta}). 
Recall  ([DeGa], III ¤6) that the scheme structure of
 $\G$ is given by the local structure theorem of Cartier-Dieudonn\'e
 $$k[\G] = k[\G_0]\otimes_kk[\G_{\red}]=
(k[X_1,\cdots,X_m]/(X_1^{p^{r_1}},\cdots,X_m^{p^{r_m}}))
 \times k[\G_{\red}]$$ 
Express the matrix of the multiplication in $k[\G]$
  by $\mu^* (a)$ on the basis $X_1^{i_1}\cdots X_m^{i_m} \otimes e_j, 
\,(0 \le i_\alpha < p^{r_\alpha})$,  the $ e_j's$ being a basis of 
$k[\G_{\red}]$.
This matrix is a sum of upper triangular blocks. Hence
 the characteristic polynomial
 of $a\in A$ is given by:  
\begin{equation}
\label{jolie}
\chi_{\G,a}(T) = \chi_{\G_{\red},a}(T)^{p^r}
\end{equation}
Assertion i. of lemma \ref{propref} follows from (\ref{jolie})
in the same way that in the constant group case.
Remark  that this result
 is valid not only for $\G$ over the base field $k$,  but also after any
 base change $k \rightarrow R$, where $R$ is a $k$-algebra.\\

Consider at last the general case: assume that 
$\G$ is a $k$-group scheme 
for $k$ any field (not perfect).
 Recall (\cite{DeGa} II ¤4) that the
 connected component of the identity of $\G$ commutes with a base 
field extension. Write same universal expressions relating the coefficients
of $\chi_{a_i}(T)$ with the coefficients of $\chi_{\sum_{i=1}^m a_i}(T)$ 
and also  over $\overline k$; we conclude as above to get i. when $\G$ is 
defined over a field $k$ or any $k$-algebra. 
\end{proof}
\begin{remark}
We 
 observe that the proof of 
lemma \ref{propref} does not rely on the finite flat
 group scheme structure: Let
 $H$ be a finite commutative free $R$-algebra. For any 
commutative $R$-algebra $A$,  define as above the characteristic 
polynomial $\chi_a (T)$, and the $\sigma_k (a)'s$  for $a \in
 H\otimes_R A$. As in the proof of lemma \ref{propref}, we 
want to express $s_k(x^{(1)}+\cdots +x^{(\ell)})'s$ in terms of 
the $s_j ({x^{(1)}}^{r_1}\cdots {x^{(\ell)}}^{r_{\ell}})$'s 
by means of polynomial relations
 with integer coefficients.
Assume $H$ to be defined over a domain 
$R_0$, i.e. $H = H_0 \otimes_{R_0} R$, which depends only on $H_0$.
 Since we have  $ H \otimes_R A = H_0 \otimes_{R_0} A$, we may 
restrict our discussion to $R$ a domain. 
Since the $s_j$'s commute with an arbitrary base
 change $A \rightarrow B$, we may assume that 
$A$  is a polynomial ring 
 in a finite number of variables over $R_0$. 
It suffices  now to 
get the polynomials relations over the fraction field $K$  of $R_0$.
 Then $H$ splits in  a finite  product of local algebras: $H =
 \prod_{\alpha = 1} ^m H_\alpha$. Performing a finite base field 
extension of $K$, we may assume that the residue field of each local
 factor $H_\alpha$ is equal to $K$. Now working with an explicit basis
 of $H$, the matrix of the multiplication by $a\in H\otimes_K A$ is a sum 
of blocks $M_\alpha$'s, each of the form  $M_\alpha = a_\alpha Id 
+ N_\alpha$, with $a_\alpha \in A$, and the matrix $N_\alpha$ being
 nilpotent. We get from  this  $$ \chi_a (T) = \prod_\alpha 
(T - a_\alpha)^{n_\alpha}$$ 
In this way the conclusion of the lemma \ref{propref}
extends to a finite flat groupo\"\i d action. 
\end{remark}
 \begin{remark} 
\label{rem23} 
If the base change  $R \to R'$  does not satisfy the condition ii of lemma
\ref{propref},  
then the canonical morphism $A^\G\otimes_R R' \rightarrow (A\otimes_R R')
 ^{\G\otimes R'}$  is in general neither injective nor surjective.
  With the Kleiman-Lonsted algebra, the situation is as seen above
 slightly bett
er. The morphism   $\phi:\Sigma_R^{\G}(A)\otimes _RR'\rightarrow 
\Sigma_{R'}^{\G\times R'}(A\otimes_RR')$ is under conditions i of
 lemma \ref{propref} always surjective. When $\G$ is infinitesimal, even if 
the action is free, $\Sigma_R^\G (A)$ may be different from $A^\G$. 
As an example, the field $k$ being of characteristic $p > 0$, take
 $A = k[\alpha_p]\otimes B$, where $B$ is a $k$ algebra of finite 
type, and the action is by left translation on $\G = \alpha_p$, 
i.e. $\Spec(A)$ is a trivial $ \G$-bundle. Then $A^\G = B$, and
 $\Sigma_k^\G (A) = k[B^p]$. 
\end{remark}   
It may be of interest
 to point out that the algebras $A^\G$ and 
$\Sigma_R^\G (A)$ are very near; first they have the same maximal
 spectrum (see Proposition 4.6  \cite{KlLo}). This leads to an 
extension of the Gabber lemma (\cite{KaMa} A7.2.1) to the present
 situation, in fact to a  finite flat groupo\"\i d.   
\begin{lemma} 
 \label{propGabber} 
i. The morphism  $\Spec (A^\G) \longrightarrow \Spec (\Sigma_R^\G (A))$
 is radicial, i.e. universally injective.\\ 
ii. If $R$ is a $\Z_p$-algebra, and if the order of $\G$ is $|\G|
 = p^r m$, with $(p,m) = 1$, then for any $a \in A^\G, 
\quad a^{p^r} \in \Sigma_R^\G (A)$;\\
 iii. Under the same hypothesis as in ii., for any base change $R \to R'$
 the morphism  $ \phi: A^\G\otimes_R R' \longrightarrow
 (A\otimes_R R') ^{\G\times R'}$ is radicial. 
\end{lemma} 
\begin{proof}   
 We first prove assertion iii. 
 The  argument of \cite{KaMa} A7.2.1 shows that there is no 
loss of generality in assuming  $R' = R/I$ to be the quotient
 of $R$ by an ideal $I$. Let $x \in (A\otimes_R R')^{\G\times R'}$.
 Lift $x$ to $a \in A$ such that $\mu^* (a) \equiv  1\otimes a $
 modulo $IR[\G]\otimes_R A + R[\G] \otimes_R IA$. We now show that
 $a^{p^r}$ lies in $\Im \phi$. The characteristic polynomial of $a$ 
satisfies $$\chi_a (T)\,\,\equiv \,\,(T - a )^{p^rm} \quad 
{\rm modulo}\quad  IA[T]$$ with $|\G|=p^rm$. By identifying
 coefficients of $T^{p^r}$, we obtain $$ \sigma_{p^rm - p^r} 
(a)\;\;\equiv\;\;    {p^rm \choose p^r} a^{p^r} \quad {\rm modulo}
 \quad IA$$ Since ${p^rm \choose p^r}\equiv m $ modulo $p$, 
 this coefficient is a unit in the $\Z_p$-algebra $R$ and $a^{p^r}
 \in \Im \phi $.  We next show that the $p^r$-th power of any
 $x \in \ker \phi$ is equal to zero. Let $a \in A^\G$  be a lift of $x$.
 Write $a = \sum \alpha_i a_i$ with $\alpha_i \in I$ and $a_i \in A$. 
On one hand, see that the coefficient of $T^j$ in the polynomial 
$\chi_a (T)$ is in $I^j A^\G$, so in $I A^\G$ if $j>0$. On the other hand,
 we have $\chi_a(T) = (T - a)^{p^rm}$; thus by identifying the coefficients
 of $T^{p^r}$ we get iii. This also proves   that the morphism induced 
by $\phi$ on the spectra is radicial.\hfill\break The proof of ii. is
 similar. For $a \in A^\G$, one has $$\chi_a (T) = (T - a) ^{p^rm} = 
T^{p^rm} - \sigma_1(a)T^{p^rm-1} + \cdots + (-1)^{p^rm} \sigma_{p^rm}(a)$$ 
By identifying the coefficients of $T^{p^r}$ in the r.h.s. and l.h.s.
 we get ${p^rm \choose p^r} a^{p^r} \in \Sigma_R^\G (A)$ as required.\\ 
Consider now i.  If $R$ is a $\Z_p$-algebra, then the result follows 
from ii. and iii. It suffices to prove that $\Spec A^\G \rightarrow 
\Spec \Sigma_R^\G (A)$ is injective. The argument of [DeGa] III,§2,3 
 (used in the proof that $\Spec A^\G$ is the topological quotient of 
$\Spec A$ by the $\G$-action)  shows without any change the injectivity
 as required. 
\end{proof} 
\begin{remark}  
If the action of the constant group $G$ is free, then we have
 $\Sigma_R^G (A) = A^G$, and the formation of $A^G$ commutes
 with an arbitrary base change.  To see this, after a faithfully 
flat base change, we may assume that $\Spec A $ is a trivial $G$-torsor,  
that is $A = R[G]\otimes_R B$ and $B = A^G$. Since $\Sigma_R^G (A)$ 
 is a surjective image of $\Sigma_R^G (R[G]) \otimes_R B = R \otimes_R
 B = B$, we have  $\Sigma_R^G (A)=B$.
 \end{remark} 
 \begin{remark}
 The group $G$ being constant, assume the action of $G$ 
is free on a dense open set. After Lemma \ref{propGabber},
 the morphism $\varphi :\Spec (A^G) \rightarrow \Spec 
(\Sigma_R^G (A))$ is radicial and is an isomorphism over 
a dense open set. If $A$ is a domain, $\varphi$ is birational.
 \end{remark} 
\subsection{Reductions}
 \label{propriete} In this paragraph, we show two properties
 of reduction of the Kleiman-Lonsted algebra: reduction to a quotient
 group scheme (Lemma \ref{dist}) and to an induced action by a group 
scheme (Lemma \ref{AInd}). We assume the conditions of lemma 2.2 fulfilled,
 i.e. either $G$ is constant or $\G$ is a $k$-group scheme and all algebras
 are over $k$.
 \begin{lemma} 
\label{dist} Assume the subgroup 
 $\H \triangleleft \G$ acts trivially on $X = \Spec A$.\\ 
i. There is an induced action of quotient group scheme $\G/\H$ 
on $X$ such that $X/\G = X/(\G/\H)$, i.e. $A^\G = A^{(\G/\H)}$.\\ 
ii. If the order $m = |\H|$ is invertible in $A$, 
then we also have $\Sigma_R^\G (A) = \Sigma_R^{\G/\H} (A)$.
 \end{lemma} 
\begin{proof}
 Assertion i. simply is the translation of the universal property of quotient scheme ([DeGa] III, §2, 6).  For assertion ii., we first establish the following relation between the characteristic polynomials of $a\in A$: 
\begin{equation} 
\label{relNm} 
\chi_{a}^\G (T) \quad = \quad \chi_{a}^{\G/\H} (T)^m  
\end{equation}
 Since the action of $\G$ factorizes through an action of $\G/\H$, 
the coaction morphism $\mu^*: A \rightarrow R[\G] \otimes_R A$ factorizes
 through the subring $R[\G/\H] \otimes_R A$: $$\mu^*:  A \longrightarrow
 R[\G/\H]\otimes_R A \,\subset \, R[\G]\otimes_R A$$ By definition $R[\G]$
 is free of rank $m$ over $R[\G/\H]$. Then taking the norm of $\mu^* (a)$,
 we get the equality $N_\G(a) = N_{\G/\H} (a)^m$. We apply this to the 
action of $\G$ on $\Spec (A[T])$, then we get the expected identity
 (\ref{relNm}).\\  Expanding the two terms of (\ref{relNm}),  we get
 the inclusion $\Sigma_R^\G(A) \subset \Sigma_R^{\G/\H} (A)$.\\  
\noindent To etablish the opposite inclusion, denote $n=|\G|$, $\ell = n/m = |\G/\H|$ and $$\chi_a^\G(T)=\sum_{j=0}^n\alpha_j^{\G}(a)T^j,\;\;\; 
\chi_a^{\G/\H}(T)=\sum_{j=0}^\ell\alpha_j^{\G/\H}(a)T^j$$ 
 After (\ref{relNm}),  $$\forall \;0\leq j\leq n,\;\;\alpha_j^\G(a)= 
\sum_{\begin{array}{cc}\scriptstyle{i_1+\cdots+i_m=j}\\ 
 \scriptstyle{ 0\leq i_1,\ldots,i_m\leq  \ell}\end{array}}
\alpha_{i_1}^{\G/\H}(a)\cdots\alpha_{i_m}^{\G/\H}(a)$$ 
In order to prove $\Sigma_R^{\G/\H}(A)\subset\Sigma_R^\G(A)$,
 it suffices to show  $$\forall \;0\leq u\leq \ell,\;\;
\alpha_u^{\G/\H}(a)\in\Sigma_R^\G(A)$$ For this, we make a
 decreasing induction on $u$. First  
$\alpha_\ell^{\G/\H}(a)=1\in\Sigma_R^\G(A)$. Assume 
 $\forall\; u<i\leq \ell,\;\;\alpha_i^{\G/\H}(a)\in\Sigma_R^\G(A)$.
 Denote $j=\ell(m-1)+u$ and consider the element of $\Sigma_R^\G(A)$
 $$\alpha_j^\G(a)=\sum_{\begin{array}{cc}\scriptstyle{i_1+\cdots+i_m=j}\\ 
 \scriptstyle{ 0\leq i_1,\ldots,i_m\leq   \ell}\end{array}}
\alpha_{i_1}^{G/H}(a)\cdots\alpha_{i_m}^{G/H}(a)$$ Since 
$j=i_1+\cdots+i_m=\ell(m-1)+u$, $0\leq i_1,\ldots,i_m\leq \ell$
 and $0\leq u\leq \ell$, the indices $i_q$, $1\leq q\leq m$ of the
 previous sum are  at least $u$. If one of them is equal to $u$, then
 all others are equal to $\ell$: $$\alpha_j^\G(a)=
m\alpha_u^{\G/\H}(a)(\alpha_l^{\G/\H}(a))^{m-1}+ \sum_{\begin{array}{cc}\scriptstyle{i_1+\cdots+i_m=j}\\ 
 \scriptstyle{ u< i_q\leq   l}\end{array}}
\alpha_{i_1}^{\G/\H}(a)\cdots\alpha_{i_m}^{\G/\H}(a)$$
 By induction, $i\in]u,\ell],\;\;\alpha_i^{\G/\H}(a)\in\Sigma_R^\G(A)$.
 Recall moreover that  $\alpha_\ell^{\G/\H}(a)=1$ and that $m$ is invertible
 in $A$. Hence $\forall u\leq i\leq \ell,\;\;\;\alpha_i^{\G/\H}(a)
\in\Sigma_R^\G(A)$. At last $\Sigma_R^\G(A)=\Sigma_R^{\G/\H}(A)$.
 \end{proof}   
  Assume that $C$ is a finite type $R$-algebra  with an action of a 
flat $R$-subgroup scheme $\H$ of $\G$.  The induced $R$-algebra of $C$
 from $\H$ to $\G$  is by definition (\cite{Ja} \S 2.12) : $$A = 
\Ind_\H^\G (C)=(R[\G]\otimes C)^\H$$ where the left action of
 $\H$ is diagonal:  $h(g,x) = (gh^{-1} , hx)$ at the points level.
 In other words, $\Spec A$ is the total space of the associated 
 $\H$-bundle $\G\times_\H  \Spec C \rightarrow \G/\H$. Note that 
the fiber over the origin $\overline{\H} \in \G/\H$ is $\Spec C$, 
 meaning there is a surjection $\phi: A \rightarrow C$,
 $\phi = e_\G\otimes 1_C$, with $e_\G$  the evaluation map. 
 \begin{lemma}
 \label{AInd} 
 i. The morphism $\phi$ induces an isomorphism $A^\G \cong C^\H$.\\
 ii. If the homogeneous space  $\G/\H$ is \'etale over $R$, then   $\Sigma_R^\G (A) \stackrel \phi\longrightarrow \Sigma_R^\H (C)$.
 \end{lemma}
 \begin{proof} i. The first assertion comes from the bundle 
interpretation of the induced $R$-algebra : the group scheme 
 $\G\times\H$ acts on $\G\times \Spec C$ via 
$$(g,h)(\sigma,c)=(g\sigma h^{-1},hc)$$ on the functor 
of points. The actions of $\G$ and $\H$  commute in 
$\G\times \H$ and $\G$ acts on $\G\times\Spec C/\H$. 
Taking the quotient by $\G$ we obtain 
$$(\G\times_\H\Spec C)/\G\simeq (\G\times \Spec C) /(\G\times \H)\simeq
 (\G\times\Spec C/\G)/\H\simeq \Spec C/\H$$ Hence $A^\G\simeq C^\H$.\\
 ii. By hypothesis $\G/\H$ is \'etale over $R$. 
After an \'etale base change,   we may  assume that $\G/\H$ splits over 
$R$.   The $R$-algebra $A$ splits into $\vert \G/\H\vert$ 
 factors as an induced module usually does, $A = \prod_{\alpha \in \G/\H}
 A_\alpha$. We can see $A_\alpha \cong C$ as the ring of functions
 with support in the $\alpha$-coset in $G$,  and $A_{\overline \H} = C$
 the fiber over the origin $\overline{\H}\in \G/\H$. From this picture,
 it is clear that if  $f\in C=A_{\overline{\H}}$, the characteristic
 polynomials satisfy  \begin{equation} \label{equaj} \chi_\G (f) = 
\prod_\alpha \alpha (\chi_\H (f)) \end{equation} Let us  denote by 
$\sigma_k^\H$ (resp $\sigma_k^\G$)  the corresponding coefficients. 
After (\ref{equaj}), $\sigma_k^\G(f) = \sigma_k^\H(f)$ for all
 $k<m=\vert  \H \vert$ and $f\in C=A_{\overline\H}$. This gives 
the inclusion $\Sigma_R^\H (C) \subset \Sigma_R^\G (A)$.  To get the
 equality, we now compute the symmetric functions $\sigma_k (f)$ for
 any $k$ and any $f\in A$.  Write $f$ as a sum $f=\sum_{\alpha\in\G/\H}
 f_\alpha \in\prod_{\alpha \in \G/\H} A_\alpha $.   As seen in the proof 
of lemma 2.2,  such a symmetric expression of $f$ is a polynomial 
expression of  symmetric functions of ``pure''  elements (belonging 
to a factor  $A_\alpha$).  For a pure element, the formula (\ref{equaj})
 gives the answer.  Hence we have the required equality.
 \end{proof} 
\noindent
 For $G$ a constant group, Lemma \ref{AInd} is well known
 (see  Raynaud (\cite{Ra4} X)).
  \section{Base change theorem} 
\subsection{Base Change} \label{critere} We first recall an elementary version of the exchange lemma (Proposition 7.7.10 \cite{Gr2}).
 \begin{lemma}
 \label{lemou} Let $R\rightarrow A$ be a noetherian local ring
 homomorphism  and $k$ be the residue field of $R$. Let $L\stackrel
 \alpha\longrightarrow P$ be a homomorphism of finite type  $R$ flat 
 $A$-modules. If the morphism
 $$\ker \alpha\otimes_R k\stackrel \phi\longrightarrow \ker
 (\alpha\otimes 1_k)$$ 
is surjective then it is bijective and the formation 
of $\ker \alpha$ is compatible with the base change: for any $R$-algebra
 $R'$, $$\ker \alpha\otimes_R R'\stackrel{\sim}\longrightarrow 
  \ker(\alpha\otimes 1_{R'})$$ \end{lemma} \noindent 
The exchange lemma \ref{lemou} allows us to prove:
 \begin{prop} 
\label{thcht}
 Let $\pi:X\rightarrow S$ be a finite type flat morphism  between
 locally noetherian schemes. Let $\G$ be a finite flat $S$-group 
scheme admissibly acting on $X$. Let $\{U_i=\Spec A_i\}_{i\in I}$ 
be a $\G$-invariant  affine cover of $X$ such that $\pi(U_i)$ is 
contained in an affine  open $S$-subscheme $V_i=\Spec R_i$ for any $i$.
 Assume that for any morphism $R_i\rightarrow \Omega$ ($\Omega$ an 
algebraically closed field) the morphism $$A_i^{\G}\otimes_{R_i}
\Omega\rightarrow (A_i \otimes_{R_i}\Omega)^{{\G}\times \Omega}$$
 is surjective. Under those conditions the morphism $X\rightarrow 
X/{\G}$ commutes with any base change. 
\end{prop} 
\begin{proof}
 Let us restrict to a $\G$-invariant open affine $U=\Spec A$ with
  $\pi(U)\subset \Spec R$. The ring $R$ is noetherian. Since 
$\pi:X\rightarrow S$ is a finite type  morphism, $A$ is a finite
 type $R$-algebra.\\ Assume $R$ to be local. Let $k$ denote its
 residue field and $\Omega$ an  algebraic  closure of $k$. 
Since $A^{\G}\otimes_R\Omega\rightarrow(A\otimes_R
 \Omega)^{\G\times\Omega}$ is surjective, $A^{\G}\otimes_Rk
\rightarrow(A\otimes_R k)^{\G\times k}$ is surjective.
 Consider the finite type homomorphism of $A$-flat $R$-module
 $$\alpha:A\rightarrow A\otimes_RR[{\G}],\;\;\; a\mapsto
 (\mu^*(a)-a\otimes1)$$ By definition $\ker\alpha=A^{\G}$.
 The hypotheses of the exchange lemma  \ref{lemou} are  satisfied. 
Hence $A^{\G}\otimes_RR'\simeq(A\otimes_R R')^{\G\times R'}$ for any
 $R$-algebra $R'$.\\ 
If the noetherian ring $R$ is not local, the same proof at any  
localization of $R$ allows us to prove proposition \ref{thcht}.
 \end{proof} 
\begin{remark} 
In proposition \ref{thcht}, it is difficult to ensure the  surjectivity 
 of morphisms $A_i^{\G}\otimes_{R_i}\Omega\rightarrow (A_i 
\otimes_{R_i}\Omega)^{{\G}\times \Omega}$ (see remark \ref{rem23}).
  \end{remark}  
The quotient  $X/\G$ is said to be  co-generated by the
 $\G$-symmetric functions on $S$ (\cite{KlLo}) if there exists a
 $\G$-invariant open affine cover $\{U_i=\Spec A_i\}_{i\in I}$  
of $X$ such that for any $i\in I$,   $A_i^{\G}=\Sigma_{R_i}^{\G}(A_i)$
  and $\pi(U_i)$ is contained in an open affine subscheme 
 $V_i=\Spec R_i$ of $S$. It is not difficult to extend this 
definition and proposition \ref{thmcht} to a finite flat groupo\"\i d.
 \begin{prop} 
\label{thmcht} Let $X\rightarrow S$ be a flat morphism and let
 $\G$ be a finite flat group scheme admissibly acting on $X$. 
Assume that 
for any geometric point $s$ of $S$, the quotient $X_s/\G_s$ 
of the geometric fiber $X_s$ is cogenerated by the 
$\G_s$-symmetric functions. Then $X/\G$ is cogenerated by 
the  $\G$-symmetric functions  on $S$ and the quotient map
 $X\rightarrow X/\G$ commutes with any base change. 
\end{prop} 
\begin{proof} 
In order to obtain proposition \ref{thmcht},  it suffices to remark 
that the hypotheses of the proposition \ref{thcht} are satisfied: 
let $s$ be a geometric point of $S$ corresponding to the algebraically
 closed field $\Omega$ : $\bar{s}:k(s)\rightarrow \Omega$. By hypothesis,
 the fiber $X_s$ is cogenerated by the $\G_s$-symmetric functions.
 Let $\{U_i=\Spec A_i\}_{i\in I}$ be a $\G$-invariant open affine 
cover of $X$ such that for any $i\in I$, $\pi(U_i)$ is contained in an
 open affine subscheme $V_i=\Spec R_i$ of $S$ and such that 
$(A_i\otimes_{R_i}\Omega)^{\G_s}=\Sigma_{\Omega}^{\G_s}(A_i)$.
 The surjectivity of the base change morphism for the invariant 
algebra  $$A_i^\G\otimes_{R_i}\Omega\rightarrow(A_i\otimes_{R_i}
\Omega)^{\G_s}$$ comes from the surjectivity of the base change
 morphism for the  Kleiman-Lonsted algebra $\Sigma_{R_i}^\G(A_i)$ 
(Lemma \ref{propref}) and from the obvious inclusion 
$\Sigma_{R_i}^\G(A_i)\subset A_i^G$. Proposition \ref{thcht} 
then proves the announced result. 
\end{proof} 
  After proposition 
\ref{thmcht}, to ensure that the quotient map commutes with an 
arbitrary base change, it is sufficient to test the cogeneracy hypothesis
 along the geometric fibers of $X \to S$. We establish two versions of 
the base change theorem: first,  for $X$ an affine smooth curve with 
an action of a finite flat group scheme (\S \ref{localversion}, local
 version) and second, for $X$ a semistable $S$-curve with an action 
 of a finite group (i.e. constant) of $S$-automorphisms  (\S 
\ref{semistable}, semistable version). 

\subsection{Local version of base change theorem} 
\label{localversion}  
  Let us begin with a local result about the action of
 finite group schemes on discrete valuation rings: 
\begin{lemma}
 \label{leminvsigma} Let $\G$  be a finite group scheme over  
the perfect field k acting  on $A = k[[t]]$, the action being 
free at the generic point.  Then $A^\G = k[[N_\G(t)]]$.
 \end{lemma} 
\begin{proof}
Let us denote the co-action of $\G$ by $\mu^*:k[[t]]\rightarrow k[\G]\otimes
k[[t]]$. The group scheme
 $\G$ is 
the semi-direct product  of its \'etale part 
and its local part of order $p^r$ (\cite{Ta}).    
The action of $\G$ being   free at 
the generic point  of $\Spec (k[[t]])$, we know the  fraction
 field  of the ring $A^\G$ has index $n = \vert \G\vert$ in the field $k((t))$ (\cite{Mu} AV). The computation made in lemma 2.2, tells us that the norm map is 
 $$ N_{\G} (f) = N_{\G_{red}} (f) ^{p^r}$$ 
 From this we see that the norm of $\mu^*(t)$
 has valuation exactly  equal to $n$. Thus $k[[t]]$ has rank $n$ over
 $k[[N_\G(t)]]$. This yields lemma  \ref{leminvsigma}.
 \end{proof}
\begin{remark}
 It should be note that the 
conclusion of lemma \ref{leminvsigma} is false in general if we have more
 than one variable. In fact for the action of the group scheme
 $\G=\alpha_p$ on $k[[x,y]]$ defined by the $p$-nilpotent vector
 field $\partial/\partial x + \partial/\partial y$, we get 
$A^\G = k[[x-y,x^p,y^p]]$, and $\Sigma_k^\G (A) = k[[x^p,y^p]]$.
\end{remark}
 \begin{prop} 
\label{localprop}
 Let $k$ be a  perfect field, and let $S=\Spec k$. 
Let $X=\Spec A$ be an affine smooth $S$-curve and $\G$ be a  finite 
flat group scheme effectively  acting on $X$, the action being free
 at the generic point of $X$.   Then $X\rightarrow X/\G$ commutes with
 any base change.
 \end{prop} 
\begin{proof} We already know that $A$ and $A^\G$ are finite module 
over $\Sigma_R^\G(A)$. Since $\Sigma_R^\G(A)=\Sigma_{\Sigma_R^\G(A)}^\G(A)$,
 we may assume without loss of generality that  $R=\Sigma_R^\G(A)$.  
Then $A$ is  now finite over  $R$ and  $R=\Sigma_R^\G(A)\subset A^\G\subset
 A$.  We must prove $R=A^\G$. The latter equality is a local property 
  on $R$, so after localization at a given  maximal ideal of $R$, we may
 assume $R$ to be local  (see lemma \ref{propref}). Then $A$ is a semi 
local ring. The hypothesis about the residue fields of $A$  implies that
 the residue field of $R$ is perfect.  The next reduction we wish to perform
 is  to pass to the completion $\hat{R}$  of $R$, which leads to a
 faithfully flat base change.  The completion $\hat{A}=A\otimes_R 
\hat R$ of $A$ satisfies $$R=\Sigma_R^\G(A)=A^\G  \mbox{ if and only if }
 \hat{R}=\Sigma_{\hat{R}}^\G(A\times_R\hat{R})=(A\times_R\hat{R})^\G$$ 
One can  now assume $R$ to be a complete local ring. The local ring 
$A^\G$ is complete and $A$  decomposes as the product of its local 
factors  indexed by the maximal ideals of the semi local ring $A$. 
 Since $\G$ acts transitively on these factors, we have $A = 
\Ind_\H^\G (A_0)$ where  $A_0$ is one of these local rings and $\H$ the
 stabilizer. We may finally confine our analysis to the case  where $A$ 
is a  complete discrete valuation ring, with perfect residue field 
(Lemma \ref{AInd}) and $\H$-action.  There exists an \'etale finite local
 extension $R \rightarrow R'$ such that  the residue fields of $A$ and 
$R'$ co\"\i ncide.  We perform the base change $R \rightarrow R'$ which saves 
the latter hypotheses: $A\otimes_R R'$ is again a product of complete 
discrete valuation rings.  We may now suppose that $R$ and $A$ have 
the same residue fields.  Now we have $A=k[[t]]$;  applying lemma 
\ref{leminvsigma} yields proposition \ref{localprop}: 
$A^\H=\Sigma_R^\H(A)=R$. \end{proof} 
 \begin{remark} 
The proof of proposition  \ref{localprop} also shows that
 $A$ is finite and flat of rank $|\G|$ over $A^\G$  and 
$A^\G$ is a co-generated Dedekind ring: $\Sigma_R^\G(A) = 
A^\G$ (i.e. a smooth curve). 
 \end{remark} 
 \subsection{Quotient of semistable curve} \label{semistable} Let us begin by describing the passage to algebra of invariants
for the  local ring at a double point (see \cite{Ra2} Appendice).\\
 Let $k$ be an algebraically closed field. Let $A={k[[x,y]]\over (xy)}$   be the completion of the local ring at a double point. In other words, $A$ is isomorphic to $$A\simeq\{(f(x),g(y))\in k[[x]]\times k[[y]],\;\; f(0)=g(0)\}$$ The group $\Gamma=\Aut_kA$ of $A$-automorphisms is the semi-direct product of $\Gamma_0$ by $(1,\tau_{x,y})$ with  $\Gamma_0=\{(\sigma,\sigma')\in \Aut k[[x]]\times\Aut k[[y]]\}$  and $\tau_{x,y}$  the $A$-automorphism exchanging the two branches: $\tau_{x,y}\in\Aut_k A$, $\tau_{x,y}(x)=y$ and $\tau_{x,y}(y)=x$.\\  Let $R$ be a noetherian ring such that $A$ is a finite type $R$-algebra. Let $G$ be a finite subgroup of $\Aut_RA$. Denote $G_0$ the maximal subgroup $G$  fixing the branches of $A$. Let $pr_x$ (resp. $pr_y$) be the projection of $A\simeq\{(f(x),g(y))\in k[[x]]\times k[[y]],\;\; f(0)=g(0)\}$ on $k[[x]]$ (resp. on $k[[y]]$). Denote by $\Delta_x$ the kernel of the projection $H$ of $G_0$ on $\Aut k[[x]]$ and by $\Delta_y$ the kernel of the projection $K$ of $G_0$ on $\Aut k[[y]]$.   \begin{lemma} \label{leminv} If the order $|\Delta_x|$ et $|\Delta_y|$ are invertible in $A$, then $A^G$ is co-generated by the $G$-symmetric functions, $A^G=\Sigma_R^G(A)$.  \end{lemma} \begin{proof}
By definition
$A^{G_0}\simeq\{(f(x),g(y))\in k[[x]]^{H}\times k[[y]]^{K},\;\;
f(0)=g(0)\}$. Hence
$$ A^{G_0}={k[[u,v]]\over (uv)}\mbox{ with } u=\Nm_Hx, v=\Nm_K y$$
The hypotheses of lemma \ref{dist} are satisfied
for the normal subgroups $\Delta_x$ and  $\Delta_y$ of
$G$ respectively acting on $pr_x(A)$ and $pr_y(A)$. Hence
$\Sigma_R^{H}(pr_x(A))= \Sigma_R^{G_0}(pr_x(A))\subset\Sigma_R^{G_0}(A)$ and 
$\Sigma_R^K(pr_y(A))=\Sigma_R^{G_0}(pr_y(A))\subset\Sigma_R^{G_0}(A)$. In 
particular
$u=\Nm_H(x)\in\Sigma_R^{G_0}(A)$ and
$v=\Nm_K(y)\in\Sigma_R^{G_0}(A)$. Then $A^{G_0}=\Sigma_R^{G_0}(A)$.\\
Hence if $G$ fixes the branches of $A$, 
proposition \ref{leminv} is proved for $G=G_0$.\\

Assume now that $G$ does not fix the branches of $A$. 
Let $\psi\in G-G_0$. 
The image $\bar{\psi}$ of $\psi$ in $G/G_0$ is an involution of
$A^{G_0}$
which exchanges the branches ; then $\bar{\psi}$ can be written as 
$\bar{\psi}=\tau_{u,v}\circ(p(u)+p^{-1}(v))$ with
$p\in\Aut k[[u]]$. An element $f(u)+g(v)\in A^{G_0}$ is
$\bar{\psi}$-invariant if and only if $g(v)=pf(v)=\bar{\psi}(f(u))$. Then
$$k[[u]]\stackrel\sim\longrightarrow A^G,\;\; u\mapsto u+\psi(u)$$
In order to prove $A^G=\Sigma_R^G(A)$, it
suffices to show $u+{\psi}(u)\in  \Sigma_R^G(A)$.\\
Let $\ell=|G_0|=n/2$ et $G=\{g_1,\ldots,g_n\}$. 
The $\ell$-th $G$-symmetric function at $x$ is
$$\sigma_\ell(x)=\prod_{1\leq i_1<\cdots<i_\ell\leq n}g_{i_1}(x)\cdots 
g_{i_\ell}(x)$$
Since $G=G_0\coprod\psi G_0$,
$$\sigma_\ell(x)=\prod_{g\in G_0}g(x)+\psi\Big(\prod_{g\in G_0}g(x)\Big)+
B_\ell(x)$$
The  sum  $B_\ell(x)$ of monomic polynomial of degree $\ell$ contains 
terms of the form $g(x)$ and $\psi g'(x)$ with $g,g'\in G^0$. Since
$g(x)$ is a serie in $x$ without any constant term and $\psi g'(x)$ 
is a serie in $y$ without any constant term, $g(x)\psi g'(x)=0$ in
$A$. This leads us to 
$$\sigma_k(x)=\Nm_{G_0}(x)+\psi(\Nm_{G_0}(x))=u+\psi(u)\in\Sigma_R^G(A)$$
and $A^G=\Sigma_R^G(A)$. 
\end{proof}
 Let $S$ be a locally
noetherian scheme. A $S$-curve
is a  flat, separated finite type relative dimension 1 $S$-scheme.  
A $S$-curve $X\rightarrow S$ is said to be {\it semistable} if
\begin{itemize}
\item $X$ is
proper,
\item its geometric fibers are 
reduced connected curves
where the singular points
are ordinary double points,
\end{itemize}
In this section $G$ is a
finite group of $S$-automorphisms on $X$.
The quotient by $G$ of the semistable curve $X$ is semistable
(\cite{Ra2} Appendice).  
In lemma \ref{leminv} we have described
the local action of $G$. We now prove our main result, 
the base change theorem \ref{thfinal}.\\

\begin{theorem}
\label{thfinal}
Assume $X$ to be a semistable $S$-curve and $G$ a finite group of
$S$-automorphisms of $X$. 
For any geometric fiber $X_s=X\times_{k(s)}\Omega$ and 
for any point $P$ of
$X_s$, we assume that the order of the kernels $\Delta_x$ and $\Delta_y$ (if $P$ 
is a double point) or $\Delta$ (if $P$ is regular) are
prime to the characteristic of $\Omega$.
Then the morphism  $(X/G)\times_S T\rightarrow (X\times_S T)/G$ is an isomorphism for any base change $T\rightarrow S$.  
\end{theorem} 
\begin{proof} 
Let $s$ be a geometric  $S$-point corresponding to the algebraic closed field  $\Omega$.  Let  $\{U_i=\Spec A_i\}_{i\in I}$ be a open                                       $G$-invariant  affine cover of $X$ such that $\pi(U_i)\subset V_i=\Spec R_i$.  The $A_i$'s are finite type $R_i$-algebra with perfect residue field at the maximal  ideals. We have to show that $(A_i\otimes_{R_i}\Omega)^G=\Sigma_{R_i}^G(A_i)$. Fix $i\in I$ and let us write $A=A_i\otimes_{R_i}{\Omega}$. We prove  $A^G=\Sigma_{R_i}^G(A)$.\\
 Performing the same reduction as in proposition \ref{localprop}, 
 this equality is  equivalent to ${A'}^{G'}=\Sigma_{R'}^{G'}(A')$ 
with $A'$ a local complete discrete valuation ring  finite on the 
local complete ring $R'=\Sigma_{R'}^{G'}(A')$ with  algebraically
 closed residue field. The algebra $A'$  comes from $A$ after a 
faithfully flat base change. The group $G'$ is the  stabilizer of
 $A'$. Since $A=A_i\otimes_{R_i}\Omega$ and $\Spec A_i$ is an affine
 open set of the  semistable curve $X$, the algebra $A'$  is 
$$A'=\Omega[[x]] \mbox{ or } {\Omega[[x,y]]\over (xy)}$$
 If $A'=\Omega[[x]]$ then ${A'}^{G'}=\Omega[[\Nm_{G'/\Delta}(T)]]=\Sigma_{R'}^{G'/\Delta}(A')$. After lemma \ref{dist} for the normal 
subgroup $\Delta$ of $G'$, $\Sigma_{R'}^{G'/\Delta}(A')=\Sigma_{R'}^{G'}(A')$. Hence   
${A'}^{G'}=\Sigma_{R'}^{G'}(A')$.\\ 
 If $A'= {\Omega[[x,y]]\over (xy)}$, ${A'}^{G'}=\Sigma_{R'}^{G'}(A')$ after lemma \ref{leminv}.\\
 Hence $(A_i\otimes_{R_i}\Omega)^G=\Sigma_{R_i}^G(A_i)$ for
any $i\in I$ and the quotient $X_s/G$ of the fiber $X_s$ is co-generated
by the $G$-symmetric functions. Theorem \ref{thfinal} then proceeds from
proposition \ref{thmcht}.
\end{proof}
Theorem \ref{th11} is an easy corollary of theorem
\ref{thfinal}:
\begin{corollary}
Let $X\rightarrow S$ be a semistable curve and let $G$ be a finite group
of $S$-automorphisms acting on $X$. If the action of $X$ is free on an open dense set on any geometric fiber, the morphism 
$(X/G)\times_S T\rightarrow (X\times_S T)/G$ is an isomorphism for
any base change $T\rightarrow S$. 
\end{corollary}

\section{Induction and restriction for deformation functors}
\label{dersection}
From here on $k$ is an algebraically closed field 
 and $S=\Spec k.$ Let $X$ be a semistable $S$-curve and $G$ be 
a finite group of $S$-automorphisms freely acting on a dense  open set of $X$.

Let us first adapt the definition of \cite{BeMe}
 of the functor of the $G$-equivariant deformations of 
$(X,G)$. 
Let $W(k)$ denote the Witt vectors ring of $k$.
Let ${\mathcal C}$ be the category of local artinian $W(k)$-algebras; 
the morphisms
are $W(k)$-morphisms of local rings. A deformation of $(X,G)$ to $A$ an
object of
${\mathcal C}$ is a $G$-equivariant isomorphism class of Galois 
covers $C\rightarrow C/G$ which further induces the identity on $X$, with:
\begin{itemize}
\item $C$ a semistable $\Spec A$-curve such that the fiber over 
the closed point of $\Spec A$ is $X\simeq C\otimes_A k$;
\item the action of $G$ on $X$ lifts to $C$ such that the isomorphism 
$X\simeq C\otimes_A k$ is $G$-equivariant.
\end{itemize}
This defines a covariant functor:
$$\Def_G(X):{\mathcal C}\rightarrow\Set,\;\;\; A\mapsto\{
\mbox{deformations of } (X,G) \mbox{ to } A\}$$
Let $H$ be a normal subgroup of $G$. The restriction of the action of 
$G$ on $X$ to $H$ defines the canonical morphism
$$\res: \Def_G(X)\rightarrow \Def_H(X)$$
After theorem \ref{thfinal} the induction morphism is well-defined
$$\ind: \Def_G(X)\rightarrow \Def_{G/H}(X/H)$$
To summarize
\begin{theorem}
Let $k$ be an algebraically closed field and $S=\Spec k$. 
Let $X$ be a semistable
 $S$-curve and $G$ be 
a finite group of $S$-automorphisms freely acting on 
a dense open set of $X$.
There exist induction and restriction morphisms between the deformation
functors:
$$
\Def_{H}(X)\stackrel \res\longleftarrow \Def_G (X)
\stackrel \ind \longrightarrow
\Def_{G/H} (X/H)$$
\end{theorem}
A local-global principle reduces the study of the deformation functor 
$\Def_G(X)$ to the study of local deformation functor at singular points and wildly ramified points (\cite{BeMe2}).

\begin{remark}  If $X$ is a smooth curve, we are able to determine the  differential  of the local restriction and induction morphisms. In the local case, we could study the functor  $D_G$ (resp. $D_H$, resp. $D_{G/H}$) of deformations of an injective morphism $G\rightarrow \Aut k[[T]]$ (resp. $H\rightarrow \Aut k[[T]]$, resp. $G/H\rightarrow\Aut k[[T]]^G$). The  tangent  spaces of these local deformation functors are isomorphic to cohomology groups (\cite{BeMe} Th. 2.2). Let $k[\varepsilon]$ be the ring of dual numbers and $k[[Y]]\simeq k[[T]]^G$.  We have $$D_G(k[\varepsilon])\simeq H^1(G,\Theta_T)\;\;\mbox{with}\;\;  \Theta_T=\{h(T){d\over dT}, h(t)\in k[[T]]\}$$ $$D_H(k[\varepsilon])\simeq H^1(H,\Theta_T),\;\;\;\; D_{G/H}(k[\varepsilon])\simeq H^1(G/H,\Theta_Y)$$ The differential $\phi_{\res}$ of the restriction morphism canonically co\"\i ncides with the restriction map between cohomology groups: $$\phi_{\res}:  H^1(G,\Theta_T)\stackrel\res\longrightarrow  H^1(H,\Theta_T)$$ For the induction case, denote the inflation map by   $\inf:H^1(G/H,\Theta_T^H)\rightarrow H^1(G,\Theta_T)$. We may show that  the differential $\phi_{\ind}$ of the induction morphism satisfies that:
$$\phi_{\ind}\circ\inf :H^1(G/H,\Theta_T^H)\rightarrow H^1(G/H,\Theta_Y)$$ 
is the natural morphism arising from the inclusion  $\Theta_T^H\subset 
\Theta_Y$.   We are  also able to made an analysis of the obstruction space  
for the $G$-deformation functor in the spirit of the Hochschild-Serre spectral 
sequence (paper in preparation). Beside  the obstructions associated to $H$
 and $G/H$, there are in general mixed obstructions (see an example in (\cite{CoKa}). 
\end{remark}  
\section {Appendix} 
We want to explain briefly for convenience of the reader how to
 get lemma \ref{result}
 used in the proof of lemma \ref{propref}. Note
 that in (\cite{KlLo}) the relations are left as an exercise.
Let $q,n\in\N^*$ and let
$x^{(\ell)}$, $1\leq \ell\leq q$ be $q$ sets of $n$ variables.
Define the partial 
polarization of the elementary symmetric functions 
\begin{equation}
\label{partial}
s_{\alpha_1,\ldots,\alpha_q} (x^{(1)},\ldots,x^{(q)}) = 
\sum (x^{(1)}_{i_1^{(1)}} \cdots x^{(1)}_{i_{\alpha_1}^{(1)}}) \cdots 
(x^{(q)}_{i_1^{(q)}} \cdots x^{(q)}_{i_{\alpha_q}^{(q)}})
\end{equation}
 where the sum is over all disjoints sequences
 $i^{(1)}_1 < \cdots < i^{(1)}_{\alpha_1} ; \ldots ; 
i^{(q)}_1 < \cdots 
< i^{(q)}_{\alpha_q}$, with $\alpha_1 + \cdots + \alpha_q \le n$. 
\begin{lemma} 
\label{result}
The partial polarization of the symmetric functions
$s_{\alpha_1,\ldots,\alpha_q} (x^{(1)},\ldots,x^{(q)})$ may be expressed as 
a polynomial with integer coefficients
in the  arguments $$s_j ( (x^{(1)})^{r_1} \cdots (x^{(q)})^{r_q})$$
 where the notation $((x^{(1)})^{r_1} \cdots (x^{(q)})^{r_q})$ 
means $((x^{(1)}_1)^{r_1}
 \cdots (x^{(q)}_1)^{r_q}, \ldots, (x^{(1)}_n)^{r_1}\cdots (x^{(q)}_n)^{r_q})$.
\end{lemma}
\begin{proof}
We work by induction on 
$\vert \alpha \vert = \sum_{i=1}^q \alpha_i$. 
Develop the product of elementary symmetric functions 
$$\prod_j
 s_{\alpha_j} (x^{(j)})= s_{\alpha_1,\ldots,
\alpha_q} (x^{(1)},\ldots,x^{(q)}) + \sum_{\beta_1,\ldots,\beta_q}
 s_{\beta_1,\ldots,\beta_q}( m_1,\ldots, m_q)$$ 
where the sum is 
over index  $\beta$ with  $\vert \beta \vert < \vert \alpha \vert$, 
and the $m_i's$ are suitable monomials in the formal variables $x^{(1)},
\cdots,x^{(q)}$. Inverting these triangular relations, we get the expected
 result.
\end{proof}
   
 \bibliographystyle{amsplain} 
 
 \end{document}